\input amstex
\documentstyle{amsppt}
\NoBlackBoxes
\document
\define\m{\medskip}
\input scrload
\define \C #1{{\scr #1}}
\define\In{\operatorname{Indet}}
\define\Prin{\operatorname{Prin}}

\define\Int{\operatorname{Int}}
\define\DR{\operatorname{DR}}
\define\IM{\operatorname{Im}}
\define \T#1{\widetilde {#1}}
\define\p{\m {\it Proof.} }
\define \my{\parindent 0cm \m \bf}
\define\eps{\varepsilon}
\topmatter
\title On Thom spaces, Massey products and non-formal symplectic manifolds\endtitle
\author Yuli Rudyak and Aleksy Tralle\endauthor
\address FB6/Mathematik, Universit\"at Siegen, 57068 Siegen, Germany\endaddress
\vskip6pt
\email rudyak\@ mathematik.uni-siegen.de,july\@ mathi.uni-heidelberg.de\endemail
\vskip6pt
\address University of Warmia and Masuria, 10561 Olsztyn, Poland
\endaddress
\vskip6pt\email tralle\@ tufi.wsp.olsztyn.pl\endemail
\vskip6pt
\rightheadtext{Thom spaces and Massey products}
\abstract
  In this work we analyze the behavior of Massey products of
closed manifolds under the blow-up construction. The results obtained
in the article are applied to the problem of constructing closed
symplectic non-formal manifolds. The proofs use Thom spaces as an
important technical  tool. This application of Thom spaces is of
conceptual interest.
\endabstract
\endtopmatter

\head Introduction \endhead

In this paper we suggest a simple general method of constructing
non-formal manifolds. In particular, we construct a large family of
non-formal symplectic manifolds. Here we detect non-formality via
essentiality of rational Massey products (we say that the Massey
product is essential if it is non-empty and does not contain the zero
element, see 1.1). In this context Thom spaces play the role of a
technical tool which allows us to construct essential Massey products
in an elegant way, Lemmas 3.4 and 3.5.
\m In greater detail, we analyze the behavior of Massey products of closed
manifolds under a blow-up
construction. The knowledge of various homotopic properties of blow-ups is
important in many areas
and, in particular, in symplectic geometry \cite{Go, M, MS, BT, TO}. In the
article we present
several results (Theorems 4.9, 4.10 and 5.2) which can be regarded as a
qualitative description of
what happens to the non-formality under the blow-up. More precisely, the
essentiality of Massey
products in the ruled submanifold yields  essential Massey products in the resulting
manifold.
\m It is clear that the blow-up procedure enables us to construct a large
class of manifolds and, on the
other hand, many classes of manifolds (complex, K\"ahler, symplectic, etc)
are invariant with respect to
blow-up. Because of this, our approach turns out to be useful in
constructing non-formal manifolds with
certain structures or other prescribed properties.
\m For example, these results have a nice application to symplectic
topology. We suggest a simple way to construct a large family of new
examples of closed symplectic non-formal (and, hence, non-K\"ahler)
manifolds, including new simply-connected examples. Note that the
problem of constructing symplectic closed non-K\"ahler manifolds (the
{\it Weinstein-Thurston problem}) was and still remains of substantial
interest in symplectic geometry \cite{BT, CFG, FMG, FLS, Go}, see [TO]
for a detailed survey. It is also worth mentioning that usually the
blow-up procedure does not change the fundamental group of ambient
manifold(s), and so our approach is applicable to manifolds with
arbitrary fundamental groups. In particular, we can produce families
of simply-connected non-formal symplectic manifolds. Note that the
construction of closed non-formal symplectic manifolds is regarded as
more difficult if one demands simple connectivity, cf. the {\it
Lupton-Oprea problem} in \cite{LO,TO}. The first simply-connected
examples have very recently appeared in \cite{BT}.

\m As another example, we mention that our method enables us to construct a
large family of algebraic varieties and K\"ahler manifolds with essential
Massey products in $\Bbb Z/p$-cohomology. (Such examples were already known,
[E]; we just emphasize that our method yields a simple construction of a
large family.)
\m It seems that this nice (and a bit surprising) application of Thom spaces
is conceptually interesting in its own right. We want also to mention that
our research was originally initiated by ideas of Gitler \cite{G}, who has
used Thom spaces in studying of homology of (complex) blow-ups.
\m We use the term ``$F$-fibration'' for any (Hurewicz) fibration with fiber $F$. Also, when
we write ``an $F$-fibration $F \to E \to B$'', it means that $E\to B$ is a
fibration with fiber $F$.

\m Throughout the paper we fix a commutative ring $R$ with the unit, and
$H^*(X)$, resp $C^*(X)$ denotes the singular cohomology $H^*(X;R)$,
resp. singular cochain complex $C^*(X;R)$ of the space $X$ unless
something else is said explicitly. Similarly, $\T H^*(X)$ denotes the
reduced singular cohomology with coefficients in $R$.
\vskip6pt
{\bf Acknowledgment.} The authors are grateful to Ivan Babenko, Hans
Baues, Sam Gitler, Peter May, James Stasheff, Manfred Stelzer and
Iskander Taimanov for valuable discussions.\vskip6pt This research was
supported by the Max Planck Institut f\"ur Mathematik and by the
Polish Research Committee(KBN).\vskip6pt We express our sincere thanks
to the anonymous referee for useful advice which improved the
presentation of this work.

\head 1. Preliminaries on Massey products \endhead

\m Throughout the section we fix a differential graded associative (in the
sequel DGA)
algebra $(A,d)$. We denote by $H(A)$ the cohomology ring of $(A,d)$, and we
always assume that $H(A)$ is commutative. Given an element $a\in A$ with
$da=0$, we denote by $[a]$ the cohomology class of $a$. So, $[a]\in H(A)$.

\m Given a homogeneous element $a\in A$, we set $\overline a=(-1)^{|a|}a$.
{\my 1.1. Definition {\rm ([K], [Ma])}.} Given homogeneous elements
$\alpha_1,\ldots, \alpha_n\in H(A)$, a {\it defining system} for
$\alpha_1, \ldots, \alpha_n$ is a family $\C X=\{x_{ij}\}$ of elements
of $A,\, 1\leq i<j\leq n+1$, with the following properties:\roster
\item $[x_{i,i+1}]=\alpha_i$ for every $i$;
\item $d x_{ij}=\dsize \sum_{r=i+1}^{j-1}\overline x_{ir}x_{rj}$ if$i+1<j<n+1$.
\endroster Consider the element $c(\C X):=\dsize \sum_{r=2}^{n}\overline
x_{1r}x_{r,n+1}$. One can prove that it is
a cocycle, so we have the class $[c(\C X)] \in H(A)$. We define
$$
\langle \alpha_1, \ldots, \alpha_n \rangle:=\left\{[c(\C X)]\bigm|\text{$\C
X$ runs over all defining systems}\right\}\i H(A).
$$
The family $\langle \alpha_1, \ldots, \alpha_n \rangle$ is called the
{\it Massey $n$-tuple product} of $\alpha_1, \ldots, \alpha_n$.

\par The {\it indeterminacy} of the Massey product
$\langle \alpha_1,\ldots, \alpha_n \rangle$ is the subset
$$
\In\langle \alpha_1, \ldots, \alpha_n \rangle:=\{x-y\bigm|x,y\in \langle
\alpha_1,\ldots, \alpha_n \rangle\}
$$
of $H(A)$.
\par We say that the Massey product $\langle \alpha_1, \ldots, \alpha_n
\rangle$ is {essential} if $\langle \alpha_1, \ldots, \alpha_n \rangle\neq
\emptyset$ and $0\notin\langle \alpha_1, \ldots, \alpha_n \rangle$.
The last condition means that there is no defining system $\C X$ with $[c(\C X)] = 0$.

{\my 1.2. Remark.} Frequently, one says about {\it non-vanishing}
Massey products instead of essential Massey products. Also, people say
that the Massey product is {\it defined}, when it is non-empty.

\m For convenience of references, we fix the following proposition. The proof
follows from the definition directly.

\proclaim{1.3. Proposition} For every morphism $f: (A,d) \to (A', d')$ of
DGA algebras and every classes $\alpha_1, \ldots, \alpha_n\in H(A)$ we have
$$
f_*\langle \alpha_1, \ldots, \alpha_n\rangle \i \langle f_*\alpha_1, \ldots,
f_* \alpha_n\rangle.$$
In particular, if $\langle \alpha_1, \ldots,
\alpha_n\rangle \neq\emptyset$ then $\langle f_*\alpha_1, \ldots, f_*
\alpha_n\rangle \neq\emptyset$, and if, in addition, $\langle f_*\alpha_1, \ldots,
f_*\alpha_n\rangle$ is essential then $\langle \alpha_1, \ldots,
\alpha_n\rangle$ is.
\qed
\endproclaim

\proclaim{1.4. Proposition} Given a DGA algebra $A$, take an element $\xi \in H(A)$
such that $\xi$ is represented by a central element $a\in A$, i.e., $ab=ba$ for every $b\in A$.
Then for every $\alpha_1, \ldots, \alpha_n\in H(A)$ and every $k$ we have
$$
\xi \langle \alpha_1, \ldots, \alpha_n \rangle \i \langle \alpha_1, \ldots,
\alpha_{k-1}, \xi \alpha_k, \alpha_{k+1}, \ldots, \alpha_n \rangle
$$
\endproclaim
\p Consider a defining system $\{x_{ij}\}$ for $\alpha_1,\ldots, \alpha_n$. We set
$$
y_{ij}=\cases ax_{ij}&{\text{ if $i\leq k$ and $j\geq k+1$,}}\\
x_{ij}& \text{otherwise}\endcases
$$
It is easy to see (using the equality $ba=ab$ for every $b\in A$) that
$\{y_{ij}\}$ is a defining system for $\alpha_1, \ldots, \alpha_{k-1},
\xi\alpha_{k}, \alpha_{k+1}, \ldots, \alpha_n$. Furthermore,
$$
\sum_{r=1}^n\overline y_{1r}y_{r,n+1}= b\sum_{r=1}^n\overline x_{1r}x_{r,n+1},
$$
and so
$$
\left[\sum_{r=1}^n\overline y_{1r}y_{r,n+1}\right]=
\xi\left[\sum_{r=1}^n\overline x_{1r}x_{r,n+1}\right].
$$
Thus,
$$
\xi\langle \alpha_1, \ldots, \alpha_n \rangle \i \langle \alpha_1, \ldots,
\alpha_{k-1}, \xi \alpha_k, \alpha_{k+1}, \ldots, \alpha_n \rangle. \quad
\qed
$$
\m  Now let us consider the special case of Massey {\it triple} products
$\langle\alpha,\beta,\gamma\rangle$. The above definition leads to the
following description. Let $a,b,c\in A$ be such that
$\alpha=[a],\beta=[b]$ and $\gamma=[c]$. Suppose that
$\alpha\beta=0=\beta\gamma$ and consider $x,y\in A$ such that
$dx=\overline a b$ and $dy=\overline b c$. Then
$\langle\alpha,\beta,\gamma\rangle$ consists of all classes of the
form $[\overline a y + \overline x c]$. Furthermore, the indeterminacy
of $\langle\alpha,\beta,\gamma\rangle$ is the set of elements of the
form $\alpha u+v\gamma$ where $u,v\in H(A)$ are arbitrary elements
with $|\alpha u|=|v\gamma|=|\alpha|+|\beta|+ |\gamma|-1$. For
convenience, we formulate the properties of the Massey triple products
as the following proposition.

\proclaim{1.5. Proposition} $\langle\alpha,\beta,\gamma\rangle \neq\emptyset$ if
and only if $\alpha\beta=0=\beta\gamma$. If $\langle
\alpha,\beta,\gamma\rangle \neq\emptyset$ then $\In \langle
\alpha,\beta,\gamma\rangle$ is the ideal $(\alpha, \gamma)$, and $\langle
\alpha,\beta,\gamma\rangle$ is a coset in $H(A)$ with respect to
$(\alpha,\gamma)$. So, $\langle\alpha,\beta,\gamma\rangle$ is essential
iff there exists $u\in \langle \alpha,\beta,\gamma\rangle$ with
$u\notin (\alpha,\gamma)$.
\qed
\endproclaim

\m Now, let $X$ be a topological space. If we wish to consider Massey
products in $H^*(X)$, we need a DGA algebra whose cohomology equals
$H^*(X)$. (It can be the singular cochain complex, or the Sullivan
model, etc, see below.) Consider a functor
$$
\Phi: \C T \to \C D
\tag{1.6}
$$
where $\C T$ is a subcategory of the category of topological spaces
and $\C D$ is the category of DGA $R$-algebras. Furthermore, we
require that $H\circ \Phi$ coincides with the (singular) cohomology
functor. Here $H$ means the cohomology functor for DGA algebras.
\m Now, given $\alpha_1, \ldots, \alpha_n\in H^*(X)$, we can consider the Massey products
$$
\langle \alpha_1, \ldots, \alpha_n \rangle_{\Phi}\i H(\Phi (X))=H^*(X)
$$
with respect to $\Phi$ (i.e. the defining system for
$\alpha_1,...,\alpha_n$ lies in $\Phi(X)$). A priori, these Massey
products depend on $\Phi$. In future we omit the subscript $\Phi$ and
write $\langle \alpha_1, \ldots,
\alpha_n \rangle$ instead of $\langle \alpha_1, \ldots, \alpha_n
\rangle_{\Phi}$, because we always say explicitly what is
$\Phi$, (i.e. which concrete $\Phi$ do we use).

\m There are three important examples of cochain functors $\Phi$. First,
given a space $X$, consider the singular cochain complex $C^*(X)$ of $X$. We equip
$C^*(X)$ with the standard (Alexander--Whitney) associative cup product
pairing. Then $(C^*(X),\delta)$ turns into a DGA algebra.
Now, since $H^*(X)=H(C^*(X), \delta)$, we use the previous construction
and define the Massey product
$$
\langle \alpha_1, \ldots, \alpha_n\rangle\i H^*(X), \quad \alpha_i\in H^*(X).
$$

\m Sullivan minimal models give us the second example. In greater detail, for every space
$X$ there is a natural commutative DGA algebra $(\C M_X,d)$ {\it over
the field of rational numbers $\Bbb Q$} which is a homotopy invariant
of $X$. Furthermore, if $X$ is a simply-connected (or more generally,
nilpotent) $CW$-space of finite type then $(\C M_X,d)$ completely
determines the rational homotopy type of $X$, see~[DGMS], [L].
\m A space $X$ is called {\it formal} if there exists a DGA-morphism$$
\rho: (\C M_X,d) \to (H^*(X;\Bbb Q),0)$$
inducing isomorphism on the cohomology level. Formality is an
important homotopy  property, since the rational homotopy type of any
nilpotent formal space can be reconstructed by some "formal" procedure
from its cohomology algebra. K\"ahler manifolds are formal
\cite{DGMS}. Examples of formal and non-formal manifolds occurring in
various geometric situations can be found in \cite{TO}.

\m Since the cochain functor $\Phi(X)=(H^*(X),0)$ yields inessential
Massey products in $H^*(X)$, we conclude that, for every formal space
$X$, all the Massey products in $H^*(X;\Bbb Q)$  with respect to the
Sullivan model are inessential, cf. \cite{DGMS, TO}. In other words, a
space $X$ is not formal if we can find an essential Massey product
(with respect to the Sullivan model) in $H^*(X;\Bbb Q)$.
\m The third example is the de Rham algebra $((\DR(X),d))$ of differential
forms on a smooth manifold $X$. We have $H(\DR(X),d)=H^*(X;\Bbb R)$, and we
can define and compute the corresponding Massey products in $H^*(X;\Bbb R)$.
\m One can prove that the Sullivan model (tensored by $\Bbb R)$ yields the
same Massey product as the de Rham complex does, see [L, Theorem III.7]. One
can also prove that the singular cochain complex with rational
coefficients yields the same Massey product as the Sullivan model does, but
the proof is more complicated: the crucial ingredients are \cite{Ma, Theorem
1.5} and \cite{BG, Prop. 3.3}. However, we do not need this fact and do not
discuss it here.

\m Since the DGA algebra of singular cochains is not commutative, Proposition
1.4 does not hold for the Massey product which comes from singular cochain
complex. However, there is the following analog of 1.4.

\proclaim{1.7. Proposition} Let the functor $\Phi$ from $(1.6)$ be
the singular cochain functor. Let $\alpha_1, \ldots \alpha_n\in
H^*(X)$ be such that  $\langle \alpha_1, \ldots,
\alpha_n
\rangle_{\Phi}\neq \emptyset$. Then for every $\xi\in H^*(X)$ and every $k$ we
have
$$
\split
\xi \langle \alpha_1, \ldots, \alpha_n \rangle \i &\pm \langle \xi\alpha_1,
\ldots, \alpha_n \rangle,\\\strut
\xi \langle \alpha_1, \ldots, \alpha_n \rangle \i &\pm \langle \alpha_1,
\ldots, \xi\alpha_n \rangle,\\\strut
\langle \alpha_1, \ldots, \alpha_{k-1}, \xi \alpha_k, \alpha_{k+1}, \ldots,
\alpha_n \rangle \,\cap &\pm \langle \alpha_1, \ldots, \alpha_{k-1},
\alpha_k, \xi\alpha_{k+1}, \ldots, \alpha_n \rangle \neq\emptyset.
\endsplit
$$
Furthermore, for Massey triple products we have:
$$
\langle \alpha, \beta, \xi\gamma\rangle \i \pm\langle \alpha, \xi\beta,
\gamma\rangle
$$
whenever $\langle \alpha, \beta, \gamma \rangle  \neq\emptyset$.
\endproclaim
Here $\pm A:=\{a|a\in A \text{ or }(-a)\in A\}$.
\p The first two inclusions and the inequality
$$
\langle \alpha_1, \ldots, \alpha_{k-1}, \xi \alpha_k, \alpha_{k+1},
\ldots, \alpha_n \rangle \,\cap\pm \langle \alpha_1, \ldots, \alpha_{k-1},
\alpha_k, \xi\alpha_{k+1}, \ldots, \alpha_n \rangle \neq\emptyset
$$
follow from \cite{K, Theorem 6}. Now, for Massey triple products we have
$$
\In \langle \alpha, \beta, \xi\gamma\rangle =(\alpha, \xi\gamma) \i (\alpha,
\gamma) = \In \langle \alpha, \xi\beta, \gamma\rangle
$$
and the desired inclusion follows from 1.5 since $\langle \alpha,
\xi\beta, \gamma\rangle \,\cap \pm \langle \alpha, \beta,
\xi \gamma\rangle \neq\emptyset$.
\qed

\head 2. $\Bbb CP^k$-fibrations\endhead

\proclaim {2.1. Proposition}Let
$$
\Bbb  CP^k @>j>> E @>p>> B
$$
be a $\Bbb  CP^k$-fibration over a path connected base, and let $\xi \in
H^2(E)$ be an element such that $j^*(\xi)\in H^2(\Bbb CP^k)$ generates the
$R$-module $H^2(\Bbb CP^k)=R$. Then the following holds:
\par {\rm (i)} every element $a\in H^*(E)$ can be represented as
$$
a=\sum_{i=0}^k \xi^ip^*(a_i),\quad a_i\in H^*(B)
$$
where the elements $a_i\in H^*(B)$ are uniquely determined by $a$;
\par {\rm (ii)} let $x, a_1, \ldots, a_m\in H^*(B)$ be such that $\xi^n
p^*x\in (p^* a_1, \ldots, p^*a_m)$ for some $n\leq k$. Then $x\in (a_1,
\ldots, a_m)$;
\par {\rm (iii)} let the functor $\Phi$ from $(1.6)$ either take the values
in the subcategory of commutative DGA algebras or is the singular
cochain complex functor. Let $\alpha, \beta, \gamma \in H^*(B)$ be
such that the Massey triple product $\langle \alpha, \beta,\gamma
\rangle$ is essential. Then the Massey triple product
$$
\langle \xi^l p^*\alpha, \xi^m p^*\beta, \xi^n p^*\gamma \rangle
$$
is essential whenever $l,m,n$ are non-negative integer numbers with
$l+m+n\leq k$.
\endproclaim
\p (i) This is the Leray--Hirsch Theorem (see e.g. [S, Theorem 15.47]).
\par (ii) We have
$$
\xi^np^*x=p^*(a_1)u_1+\cdots+p^*(a_m)u_m, \quad u_i\in H^*(E).
$$
Because of (i), every $u_i$ can uniquely be represented as
$u_i=\sum_{j=0}^k\xi^jp^*u_{ij}$ with $u_{ij}\in H^*(B)$. So,
$$
\xi^np^*x=\sum_{r=1}^mp^*a_r \sum_{j=0}^k \xi^jp^*u_{rj}=
\sum_{j=0}^k\xi^{j}\sum_rp^*(a_ru_{rj}).
$$
Now, by (i), we conclude that $x=\sum_{r}a_ru_{rn}$, i.e. $x\in
(a_1,\ldots, a_n)$.
\par (iii) First, notice that
$$
\langle \xi^l p^*\alpha, \xi^m p^*\beta,
\xi^n p^*\gamma \rangle\neq \emptyset
$$
by 1.5 and 1.3. Take $u\in \langle \alpha,\beta,\gamma\rangle$ with
$u\notin (\alpha,\gamma)$. Then, by (ii),
$\xi^{l+m+n}p^*u\notin(p^*\alpha, p^*\gamma)$, and hence
$$
\xi^{l+m+n}p^*u\notin(\xi^lp^*\alpha, \xi^np^*\gamma)\subset (p^*\alpha, p^*\gamma).
$$
So, because of 1.5, it suffices to prove that
$$
\xi^{l+m+n}p^*u\in \pm \langle \xi^l p^*\alpha, \xi^m p^*\beta, \xi^n
p^*\gamma \rangle.
$$
Now, if $\Phi(E)$ is a commutative $DGA$ algebra then, by 1.4,
$$
\xi^{l+m+n}p^*u\in \xi^{l+m+n} \langle p^*\alpha, p^*\beta, p^*\gamma \rangle
\i \langle \xi^l p^*\alpha, \xi^m p^*\beta, \xi^n p^*\gamma \rangle
$$
If $\Phi(X)=(C^*(X), \delta)$ then, by 1.7,
$$
\split
\xi^{l+m+n}p^*u\in \xi^{l+m+n} \langle p^*\alpha, p^*\beta, p^*\gamma
\rangle & \i \pm \langle \xi^l p^*\alpha, p^*\beta, \xi^{m+n} p^*\gamma\rangle\\
& \i \pm \langle \xi^l p^*\alpha, \xi^m p^*\beta, \xi^n p^*\gamma \rangle.
\quad
\qed
\endsplit
$$

\head 3. Thom spaces and Massey products\endhead

In this section the functor $\Phi$ is an arbitrary functor as in (1.6).

\m We need some preliminaries on Thom spaces of normal bundles. Standard
references are [B], [R]. Let
$M$, $X$ be two closed smooth manifolds, and let $i: M \to X$ be a smooth
embedding. Let $\nu$ be the
normal bundle of $i: M \i X$, $\dim \nu=d$, and let $T\nu$ be the Thom space
of $\nu$. We assume
that $\nu$ is orientable, choose an orientation of $\nu$ and denote by $U\in
\T H^d(T\nu)$ the Thom class of $\nu$.

\m Let $N$ be a closed tubular neighborhood of $i(M)$ in $X$. Let $V$ be the
interior of $N$, and set $\partial N = N \setminus V$. The Thom space
$T\nu$ can be identified with $X/X\setminus V = N/\partial N$. We
denote by
$$ c: X@>\text{quotient}>> X/(X\setminus V) = T\nu
\tag{3.1}
$$
the standard collapsing map.

\m Recall that, for every two pairs $(Y,A),\, (Y,B)$ of topological spaces,
there is a natural cohomology pairing
$$
H^i(Y,A) \otimes H^j(Y,B) \to H^{i+j}(Y,A\cup B)
$$
see [D]. In particular, we have the pairings
$$
\varphi: \T H^i(T\nu)\otimes H^j(M) \to \T H^{i+j}(T\nu)
$$
of the form
$$
H^i(N, \partial N) \otimes H^j(N)\to H^{i+j}(N, \partial N)
$$
and the pairing
$$
\psi:\T H^i(T\nu)\otimes H^j(X) \to \T H^{i+j}(T\nu)
$$
of the form
$$
H^i(X, X\setminus V)\otimes H^j(X)\to H^{i+j}(X, X\setminus V).
$$
It is well known and easy to see that the diagram
$$
\CD
\T H^*(T\nu) \otimes H^*(M) @>\varphi>>\T H^*(T\nu)\\@A1\otimes i^*AA @|\\
\T H^*(T\nu) \otimes H^*(X) @>\psi>> \T H^*(T\nu)\\
@V\eps \otimes 1VV @VV\eps V\\H^*(T\nu) \otimes H^*(X) @. H^*(T\nu)\\
@Vc^* \otimes 1VV @VVc^* V\\H^*(X) \otimes H^*(X) @>\Delta>> H^*(X)
\endCD
\tag{3.2}
$$
commutes; here $\Delta=\Delta_X$ is induced by the diagonal $X \to X \times
X$ and $\eps$ is the canonical inclusion $\T H^*(T\nu)\i H^*(T\nu)$.

\m As usual, for the sake of simplicity we denote each of the products $\varphi(a\otimes b)$,
$\psi(a\otimes b)$ and $\Delta(a\otimes b)$ by $ab$. Also, below we
will use the same letter for an element $x\in \T H^*(T\nu)$ and its
image $\eps(x)$ under the canonical inclusion $\T H^*(T\nu)\i
H^*(T\nu)$. For example, we can consider the Thom class $U\in
H^d(T\nu)$.

\m Finally, we recall that the Euler class $\chi=\chi(\nu)$ of $\nu$ is
defined as $\chi:=\frak z^*(U)\in H^d(M)$ where $\frak z: M \to T\nu$
is the zero section of the Thom space. Furthermore,
$$
\frak z^*(Ua)=\chi a\quad \text{for every $a\in H^*(M)$},
\tag{3.3}
$$
see e.g.~[R, Prop. V.1.27].

\proclaim{3.4. Lemma} If $\langle \alpha, \beta, \gamma\rangle
\neq\emptyset$ for some $\alpha, \beta,\gamma \in H^*(M)$ then
$\langle \chi \alpha, \chi\beta, \chi \gamma \rangle \neq\emptyset$.
If, in addition, $\langle
\chi \alpha, \chi\beta, \chi \gamma \rangle$ is essential then there are
$u,v,w\in H^*(X)$ such that the Massey product $\langle u,v, w\rangle$ is
essential.
\endproclaim
\p Clearly, $\chi \alpha \chi \beta= 0 = \chi \beta \chi \gamma$, and so
$\langle \chi\alpha, \chi\beta, \chi\gamma\rangle \neq\emptyset$. Set
$u:=c^*(U\alpha),\, v:=c^*(U\beta),\, w:=c^*(U\gamma)$ and prove that
$uv=0=vw$. Consider the diagram (where $H$ denotes $H^*$)
$$
\CD
\T H(T\nu) \otimes H(M) \otimes \T H(T\nu) \otimes H(M) @>\varphi \otimes
\varphi >> \T H(T\nu) \otimes H(T\nu) @>\Delta >>\T H(T\nu)\\ @VVTV @. @|\\
\T H(T\nu) \otimes \T H(T\nu) \otimes H(M) \otimes H(M)
@>\Delta \otimes \Delta >> \T H(T\nu) \otimes H(M) @>\varphi >>\T H(T\nu)
\endCD
$$
where $T(a\otimes b \otimes c \otimes d)=a\otimes c\otimes b\otimes d$. This
diagram commutes up to sign,  and so
$$
\split
(U\alpha)(U\beta)&=\Delta((U\alpha)\otimes(U\beta))=\Delta\circ(\varphi\otimes
\varphi)(U\otimes\alpha\otimes U\otimes\beta)\\ &=\pm \varphi(\Delta \otimes
\Delta)(U\otimes U
\otimes \alpha \otimes \beta)=\pm \varphi((UU)\otimes (\alpha\beta))=0
\endsplit
$$
since $\alpha\beta=0$. So, $uv=c^*((U\alpha)(U\beta))=0$.
\par Similarly, $vw=0$, and so $\langle u,v,w\rangle \neq\emptyset$.
Furthermore, the map
$$
M @>i>> X @>c >> T\nu
$$
coincides with the zero section $\frak z: M \to T\nu$. Now,
$$
0\notin \langle \chi \alpha, \chi\beta, \chi \gamma\rangle
=\langle \frak z^*(U\alpha), \frak z^*(U\beta), \frak z^*(U\gamma)\rangle
=\langle i^*u, i^*v, i^*w\rangle,
$$
and the result follows from 1.3.
\qed

\proclaim{3.5. Lemma} Let $\alpha, \beta\in H^*(M)$ and $w\in H^*(X)$ be
such that $\langle \alpha, \beta, i^*w\rangle \neq\emptyset$. Then $\langle
\chi \alpha, \chi\beta, i^*w \rangle \neq\emptyset$. If, in addition,
$\langle \chi \alpha, \chi\beta, i^*w \rangle$ is essential
then there are $u,v\in H^*(X)$ such that the Massey product $\langle u,
v, w\rangle$ is essential.
\endproclaim
\p Clearly, the Massey product
$\langle \chi \alpha,
\chi\beta, i^*w \rangle$ is not empty. As in the proof of 3.4, we set
$u:=c^*(U\alpha),\, v:=c^*(U\beta)$. We have proved in 3.4 that
$uv=0$. Now we prove that $vw=0$, i.e. that $\langle u,v, w\rangle
\neq\emptyset$. Consider the commutative diagram
$$
\CD
\T H^*(T\nu) \otimes H^*(M) \otimes H^*(M) @>1\otimes \Delta>> \T H^*(T\nu)
\otimes H^*(M)\\@V\varphi \otimes 1VV @VV\varphi V\\
\T H^*(T\nu)\otimes H^*(M) @>\varphi >> \T H^*(T\nu)\\
\endCD
$$
Notice that
$$
\eqalign{
\varphi(U\beta \otimes i^*w)&=\varphi\circ(\varphi\otimes 1)(U\otimes \beta
\otimes i^*w)=\varphi(1\otimes\Delta)(U\otimes \beta \otimes i^*w)\cr
&=\varphi(U\otimes \beta i^*w)=0.\cr}
\tag{3.6}
$$
Furthermore, in the diagram (3.2) we have
$$
vw=\Delta(v\otimes w)=\Delta(c^*(U\beta)\otimes w)=c^*\psi ((U\beta)\otimes
w)=c^*\varphi(U\beta\otimes i^*w)=0
,$$
the last equality follows from (3.6). Now the proof can be completed just as
the proof of 3.4.
\qed
{\my 3.7. Remark.} Certainly, in 3.5 we can consider the Massey
product $\langle i^*w, \alpha, \beta\rangle$, resp. $\langle \alpha,
i^*w,\beta, \rangle$, and get the essential Massey product $\langle w,
u, v\rangle$, resp. $\langle u, w, v, \rangle$. We leave it to the
reader to formulate and prove the corresponding parallel results.

\head 4. Applications to blow-up \endhead

In this section the functor $\Phi$ is assumed to be as in 2.1(iii).

{\my 4.1. Recollection.} Let $\zeta$ be a $(k+1)$-dimensional complex
vector bundle over a space $Y$. We assume that $\zeta$ is equipped
with a Hermitian metric and let $\Prin(\zeta)=\{P \to Y\}$ be the
corresponding principal $U(k+1)$-bundle. Recall that the {\it
projectivization} of $\zeta$ is a locally trivial $\Bbb
CP^k$-bundle
$$
\T Y\to Y
$$
where $\T Y:= P\times_{U(k+1)}\Bbb  CP^k$ and the $U(k+1)$-action on
$\Bbb CP^k$ is induced by the canonical $U(k+1)$-action on $\Bbb
C^{k+1}$.
\par Now we define a canonical complex line bundle $\lambda_{\zeta}$ over
$\T Y$ as follows. The total space $L$ of the canonical line bundle
$\eta_k$ over $\Bbb CP^k$ has the form
$$ L=\{(z,l)\in \Bbb
C^{k+1}\times \Bbb C P^k\bigm| z\in l\}.
$$
In particular, the projection $\pi: L \to \Bbb  CP^k,\, (z,l) \mapsto
l$ is an $U(k+1)$-equivariant map. We define $\lambda_{\zeta}$ to be
the induced map
$$
1\times \pi: P\times_{U(k+1)}L \to
P\times_{U(k+1)}\Bbb CP^k=\T Y.
$$
Notice that if $Y$ is the one-point space then $\T Y=\Bbb CP^k$ and
the canonical bundle $\lambda$ coincides with $\eta_k$.

\m The construction $\lambda_{\zeta}$ is natural in the following sense. Let
$\zeta$ be a complex vector bundle over $Y$, and let $f: Z \to Y$ be
an arbitrary map. Then the obvious map $I: f^*(\Prin \zeta)\to \Prin
\zeta$ yields a commutative diagram
$$
\CD
\T Z @>\T f>> \T Y\\@VVV @VVV\\Z @>f>> Y
\endCD
$$
where $\T Z \to Z$ is the projectivization of $f^*\zeta$ and $\T f$ is
induced by $I$.

\proclaim{4.2. Proposition} $\T f^*\lambda _{\zeta}$ is naturally isomorphic
to $\lambda_{f^*\zeta}$.
\qed
\endproclaim

\proclaim{4.3. Corollary} If $\Bbb  CP^k @>j>> \T Y @>>> Y$ is the
projectivization of a complex vector
bundle $\zeta$, then $j^*\lambda_{\zeta}$ is isomorphic to the canonical
line bundle $\eta_k$ over $\Bbb CP^k$.
\qed
\endproclaim

{\my 4.4. Definition.} Let $M$ and $X$ be two closed connected smooth
manifolds, and let $i: M \to X$ be a smooth embedding of codimension
$2k+2$. A {\it blow-up along $i$} is a commutative diagram of smooth
manifolds and maps
$$\CD\Bbb  CP^k @.\\@VjVV @.\\
\T M @>\T i >> \T X\\@ VpVV @VVqV\\ M @>i>> X
\endCD
\tag{4.5}
$$
such that the following holds:
\par (i) $\Bbb  CP^k @>j>>\widetilde M @ >p>> M$ is a locally trivial
bundle (with fiber $\Bbb  CP^k$) which is a projectivization of a complex
$(k+1)$-dimensional
vector bundle $\zeta$. In particular, $\T M$ is a closed connected manifold;
\par (ii) $\T X$ is a closed connected manifold, $\widetilde i:
\widetilde M \to \widetilde X$
is a smooth embedding of codimension $2$, and the line bundle
$\lambda_{\zeta}$ is isomorphic (as a
real vector bundle) to the normal bundle $\nu$ of $\T i$;
\par (iii) there is a closed tubular neighborhood $N$ of $i(M)$ such that
$\T N:= q^{-1}(N)$ is a tubular neighborhood of $\T i(\T M)$ and
$$
q|_{\T X\setminus \Int \T N}: \T X\setminus \Int \T N \to X \setminus \Int N
$$
is a diffeomorphism.

\m By 4.4(ii), the normal bundle $\nu$ of $\T i$ is isomorphic to
$\lambda_{\zeta}$, and hence $\nu$ is orientable. Take an orientation of
$\nu$ and consider the Euler class $\chi(\nu)\in H^2(\T M)$.

\proclaim{4.6. Proposition} The class $j^*\chi(\nu)$ generates the
$R$-module $H^2(\Bbb CP^2)=R$.
\endproclaim
\p  By 4.4(ii) and 4.3, $j^*\nu$ is isomorphic to $\eta_k$, and hence
$j^*\chi(\nu)=\chi(\eta_k)$.
So, the result holds for $R=\Bbb Z$, and thus it holds for arbitrary $R$.
\qed

\proclaim{4.7. Proposition} Consider a blow-up diagram $(4.5)$. If $k\ge 1$
then $q_*:\pi_1(\T X) \to\pi_1(X)$ is an isomorphism.
\endproclaim
\p Let $V$ denote the interior of $N$ and set $\partial N = N\setminus V$.
Similarly, $\T V:=q^{-1}(V)$
and $\partial \T N:=\T N \setminus \T V$. Since $k\geq 1$, the inclusion
$X\setminus V \i X$ induces
an isomorphism of fundamental groups. So, we must prove that the inclusion
$\T X\setminus \T V \i \T
X$ induces an isomorphism of fundamental groups. But this follows from the
van Kampen Theorem, since
the inclusion $\partial \T N \i \T N$ induces an epimorphism $\pi_1(\partial
\T N) \to \pi_1(\T N)$.
The last claim holds, in turn, because $\partial \T N \to \T M$ is a
locally trivial bundle with connected fiber.
\qed
{\my 4.8. Remark.} Usually the term ``blow-up'' is reserved for a
canonical procedure which, in particular, leads to a diagram like
(4.5), cf. [G], [ M], [MS]. In Definition 4.4 we have just axiomatized
certain useful (for us) properties.

\proclaim{4.9. Theorem} Consider a blow-up diagram $(4.5)$. Suppose that
$k\geq 3$ and that $M$ possesses an essential Massey triple product.
Then $\T X$ possesses an essential Massey triple product.
\endproclaim
\p Let $\langle \alpha, \beta, \gamma\rangle$ be an essential Massey
product in $M$, and let $\chi$ be the Euler class of the normal bundle
of $\T i$. Then, by 4.6 and 2.1(iii), the Massey triple product
$$
\langle \chi p^*\alpha, \chi p^*\beta, \chi p^*\gamma\rangle\subset
H^*(\T M)
$$
is essential. Thus, by 3.4, $\T X$ possesses an essential Massey triple
product.
\qed

\m Sometimes it is useful to replace the assumption $k\ge 3$ in 4.9 by a
weaker assumption $k\ge 2$. We do it as follows.
\proclaim{4.10. Theorem} Consider a blow-up diagram as in $(4.5)$. Suppose
that $k\geq 2$ and that there are
elements $\alpha,\beta \in H^*(M)$ and $w\in H^*(X)$ such that at least one
of the Massey triple
product $\langle \alpha,\beta,i^*w\rangle $, $\langle \alpha,
i^*w,\beta,\rangle $, $\langle
i^*w,\alpha,\beta,\rangle $ is essential. Then $\T X$ possesses an essential
triple
Massey product.
\endproclaim
\p We consider the case when the Massey product $\langle
\alpha,\beta,i^*w\rangle \subset H^*(M)$ is essential, all the other cases can be
considered similarly. Let $\chi$ be the Euler class of the
normal bundle of $\T i$. Then, by 4.6 and 2.1(iii), $\langle \chi p^*\alpha,\chi
p^*\beta, \T i^*q^*w\rangle \neq\emptyset$, and the Massey product
$$
\langle \chi p^*\alpha, \chi p^*\beta, \T i^*q^*w\rangle \subset H^*(\T M)
$$
is essential. Thus, by 3.5, $\T X$ possesses an essential Massey product of
the form $\langle u,v,q^*w\rangle$.
\qed

\head 5. Application to symplectic manifolds\endhead

In this section the functor $\Phi$ is assumed to be the Sullivan minimal model.

\proclaim{5.1. Theorem {\rm ([M])}} Let $M$ and $X$ be two closed symplectic
manifolds, and let
$i: M \to X$ be a symplectic embedding. Then there exists a blow-up diagram
where the map $\T i$ is
a symplectic embedding. In particular, $\T X$ is a symplectic manifold.\qed
\endproclaim
\m We define a {\it symplectic blow-up} to be a blow-up with $i$ and $\T i$
symplectic.
\m Actually,   McDuff~[M] suggested a canonical construction of symplectic
blow-up. In particular,
the bundle $\zeta$ in 4.4(i) turns out to be the  normal bundle of the
embedding $i$.
 A detailed exposition of this construction can be found in [MS], [TO].

\m This theorem enables us to use the above Theorems 4.9, 4.10 in order
to construct non-formal symplectic manifolds. To start with, we must
have at least one symplectic manifold with an essential Massey triple
product. {\my 5.2. Example {\rm(Kodaira--Thurston)}.} Let $H$ be the
Heisenberg group, i.e. the group of the $3\times 3$-matrices of the
form
$$
\alpha=\pmatrix 1 & a & b\\0 & 1 & c\\0 & 0 & 1
\endpmatrix
$$
with $a,b,c\in \Bbb R$, and let $\Gamma$ be the subgroup of $H$ with
integer entries. We set
$$
K:=(H/\Gamma) \times S^1.
$$
So, the Kodaira--Thurston manifold $K$ is a 4-dimensional nilmanifold. Its Sullivan
minimal model has the form
$$
(\Lambda(x_1,x_2,x_3,x_4),d)\quad \hbox{with}\quad dx_1=dx_2=dx_4=0,\,
dx_3=x_1x_2,
$$
where $\deg\,x_i=1$ for $i=1,2,3,4$ and the generators $x_1, x_2, x_3$
come from the Heisenberg manifold. We set $\alpha=[x_1]$, $\beta=[x_2]$, $\gamma=x_4$ and $u=[x_1x_4+x_2x_3]$. Then
$$
\split
H^1(K;\Bbb Q)&=\Bbb Q^3=\{\{\alpha,\beta,\gamma\}\},\\
H^2(K;\Bbb Q)&=\Bbb Q^4=\{\{\alpha\gamma,\beta\gamma, u, [x_1x_3]\}\},\\
H^3(K;\Bbb Q)&=\Bbb Q^3=\{\{\alpha u,\gamma u, [x_1x_3x_4]\}\}
\endsplit
$$
where $\{\{\ldots\}\}$ denotes ``the $\Bbb Q$-vector space with the basis ...''. Finally, $K$ possesses a symplectic form $\omega$, and the class
$$
u \otimes 1\in H^*(K;\Bbb Q)\otimes \Bbb R = H^*(K;\Bbb R)
$$
coincides with the de Rham cohomology class of $\omega$. See [TO] for
details.

\m It is easy to see that the Massey product $\langle \alpha,\alpha,\beta\rangle$ is essential, but we need more.

\proclaim{5.3. Proposition} The Massey products $\langle \alpha,\alpha,\beta\rangle$ and $\langle\beta,\beta,u\rangle$ are essential.
\endproclaim
\p We prove the essentiality of the second Massey product, because the proof for the first one is simpler. Clearly $\beta\beta=0$. Furthermore,
$$
\beta u=[x_2x_1x_4+x_2x_2x_3]=[x_2x_1x_4]=[-d(x_3x_4)]=0,
$$
 and hence  $\langle\beta,\beta,u\rangle\neq \emptyset$. Furthermore,
$$
x_2(x_1x_4+x_2x_3)=-d(x_3x_4),
$$
and  hence $=[-x_1x_3x_4] \in \langle\beta,\beta,u\rangle$. But $[x_1x_3x_4]\notin (\beta,u)$ because 
$$
(\beta, u)\cap H^3(K,\Bbb Q)=\{\{\alpha u,\gamma u\}\}.
$$ 
Thus, $\langle\beta,\beta,u\rangle$ is essential.
\qed
\proclaim{5.4. Theorem}Consider a symplectic embedding $i:K \to X$ and any
symplectic blow-up
$$
\CD
\Bbb  CP^k @.\\@VjVV @.\\\T K @>\T i >> \T X\\
@ VpVV @VVqV\\K @>i>> X
\endCD
$$
along $i$. If $k\geq 2$ then $\T X$ possesses an essential rational Massey
triple product. In particular, $\T X$ is not a formal space.
\endproclaim
\p For $k\geq 3$ this follows from 4.9, because $K$ possesses an essential
Massey triple product. If $k=2$, consider the symplectic form
$\omega_X$ on $X$. Then $i^*\omega_X=\omega_K$ where $\omega_K$ is the
symplectic form on $K$. Hence,
$$
u \in \IM\{i^*: H^2(X;\Bbb Q) \to H^2(K;\Bbb Q)\}.
$$
So, according to 5.3, $\langle \beta,\beta, i^*w\rangle$ is essential
for some $w\in H^2(X;\Bbb Q)$. Now, the result follows from 4.10.
\qed
\m In particular, if the initial ambient space $X$ is simply-connected, then
the space $\T X$ in 5.4 gives us an example of non-formal
simply-connected symplectic manifold.

\m It is well known that, for every $n\ge 5$, there exists a symplectic
embedding $K \hookrightarrow \Bbb CP^n$, \cite{Gr, T}. Consider a
symplectic blow-up along such an embedding, and let $\T {\Bbb CP^n}$
denote the corresponding space (in the top right corner of the diagram
(4.5)). Then 5.2 yields the following corollary.

\proclaim{5.5. Corollary {\rm (Babenko--Taimanov~[BT])}} Every space
$\T{\Bbb CP^n},\, n \ge 5$ possesses an essential Massey triple
product, and hence it is not a formal space.
\qed
\endproclaim
{\my 5.6. Remarks.} 1. Generalizing 5.2, consider the so-called Iwasawa
manifolds. Namely, we set
$$I(p,q)=(H(1,p)\times H(1,q))/\Gamma
$$
where $H(1,p)$ consists of all matrices of the form
$$
\alpha=\pmatrix
I_p & A & C \\0 & 1 & b\\0 & 0 & 1\endpmatrix
$$
It is proved in [CFG]
that $I(p,q)$ always possesses an essential rational Massey triple
product. So, we can use these manifolds for construction of families
of non-formal symplectic manifolds. In particular, if both $p$ and $q$
are even then $I(p,q)$ is symplectic, and so in this way we can also
get a large family of non-formal (and so non-K\"ahler) symplectic
manifolds.

\m 2. Certainly, the blow-up construction can be iterated, i.e. starting
from $M\i X$ and getting $\T X$, we can embed $\T X$ into another
manifold and construct a new blow-up, and so on. So, here we really
have a lot of possibilities to construct manifolds with essential
Massey products and, in particular, non-formal manifolds.

\head 6. Massey products in algebraic and K\"ahler manifolds \endhead

In this section $\Phi^*(X)=(C^*(X;\Bbb Z/p), \delta)$ where $p$ is an odd
prime number.

\m According to [DMGS], every K\"ahler manifold $M$ is formal, and so every
Massey product in $H^*(M;\Bbb Q)$ is inessential. However, there are
complex projective algebraic varieties (and in particular K\"ahler
manifolds) $S$ with essential Massey products in $H^*(S;\Bbb Z/p)$.
Indeed, Kraines~[K] proved that, for every prime $p>2$ and every $x\in
H^{2n+1}(X;\Bbb Z/p)$, we have
$$
\beta P^n(x)\in -\langle \undersetbrace p\ \text{times}\to{x,x,\ldots , x}
\rangle\subset H^{2np+2}(X;\Bbb Z/p)
$$
where $\beta: H^i(X;\Bbb Z/p) \to H^{i+1}(X;\Bbb Z/p)$ is the Bockstein
homomorphism and
$$P^n : H^i(X;\Bbb Z/p) \to H^{i+2n(p-1)}(X;\Bbb Z/p)
$$
is the reduced Steenrod power. For instance, if $H_1(X;Z)=\Bbb Z/3$
then $H^1(X;\Bbb Z/3)=\Bbb Z/3$ and $\beta(x)\ne 0$ for every $x\in
H^1(X;\Bbb Z/3),\,x\ne 0$. So, the Massey product $\langle
x,x,x\rangle$ has zero indeterminacy, and $\langle x,x,x\rangle\ne 0$
provided $x\ne 0$.

\m So, if $H_1(X;\Bbb Z)=\Bbb Z/3$, then $H^*(X;\Bbb Z/3)$ possesses an
essential Massey triple product. Certainly, one can find many examples
of algebraic varieties and K\"ahler manifolds with such homology
groups, see e.g.~[ABCKT]. Now, since the classes of K\"ahler manifolds
and algebraic varieties are invariant under the blow-up construction,
we can use Theorems 4.9 and 4.10 and construct a large class of
K\"ahler manifolds and algebraic varieties with essential Massey
triple $\Bbb Z/3$-products, including simply-connected objects as
well.

\m Perhaps, the following observation looks interesting. If we have a
K\"ahler manifold or algebraic variety $ X_0$ with $u_0\in
H^1(X_0;\Bbb Z/3)$ such that $0\notin\langle u_0, u_0,u_0\rangle$,
then we can perform a blow-up along $X_0$ (in the corresponding
category) and get a resulting object $X_1$, and if $\dim X_1-\dim X
_0\geq 6$ then there is an element $u_1\in H^3(X_1;\Bbb Z/3)$ with
$0\notin\langle u_1, u_1,u_1\rangle$. In particular, $\beta
P^1(u_1)\ne 0$. Similarly, performing an obvious induction, we can
construct $X_n$ (algebraic or K\"ahler) and an element $u_n\in
H^{2n+1}(X_n;\Bbb Z/3)$ with $\beta P^n(u)\neq 0$.

\m Can one generalize these results for modulo $p$ cohomology? One can
construct a K\"ahler manifold with the fundamental group $\Bbb Z/p$,
but we cannot go ahead because we are not able to prove an analog of
2.1(iii) for $p$-tuple Massey products. However, there is the
following analog of the above results.

\m For every $u\in H^*(X;\Bbb Z/p)$, Kraines~[K] defined a ``small" Massey
product
$$
\langle u\rangle^k \subset \langle\undersetbrace k\
\text{times}\to { u,u,\ldots , u} \rangle
$$
by considering special defining systems where $x_{ij}=x_{kl}$ whenever
$j-i=l-k$. Moreover, he proved that
$$
-\beta P^i(u)= \langle u \rangle ^p\quad  \text{ for every } u\in
H^{2i+1}(X;\Bbb Z/p).
$$
It is remarkable that $\langle u \rangle ^p$ has zero indeterminacy.
Now we can prove an analog of 2.1(iii) for $\langle u \rangle ^p$
provided $k\geq p$, i.e. to prove that $\langle \xi u \rangle ^p\ne 0$
provided $\langle u
\rangle ^p\ne 0$. (One can prove it directly, or can use \cite{K, Corollary
7}.) Thus, we can proceed (by blowing up) and construct a large family of
K\"ahler manifolds with $\langle u\rangle^p\ne 0$.

\enddocument
\end